\documentclass[a4paper]{article}

\usepackage{geometry}
\usepackage{authblk}
\usepackage{orcidlink}
\usepackage{graphicx}%
\usepackage{multirow}%
\usepackage{amsmath,amssymb,amsfonts}%
\usepackage{amsthm}%
\usepackage{mathrsfs}%
\usepackage[title]{appendix}%
\usepackage{xcolor}%
\usepackage{textcomp}%
\usepackage{manyfoot}%
\usepackage{booktabs}%
\usepackage{algorithm}%
\usepackage{algorithmicx}%
\usepackage{algpseudocode}%
\usepackage{listings}%
\usepackage{tikz}
\usetikzlibrary{arrows.meta}
\usepackage{mathtools}
\usepackage{pgfplots}
\usepackage{pgfplotstable}
\usepackage{anyfontsize}
\usepackage{subcaption}
\usepackage{placeins}
\usepackage{arydshln}
\usepackage{stmaryrd}
\usepackage{nicematrix}
\usepackage{cleveref}
\usepackage{xspace}
\usepackage{verbatim}

\pgfplotsset{width=\textwidth,
			 compat=1.9,
			 colormap={rainbow}{rgb255=(255,0,0),
                                rgb255=(255,165,0),
                                rgb255=(0,128,0),
                                rgb255=(0,0,255),
                                rgb255=(128,0,128)},
             cycle multiindex* list = {[of colormap=rainbow]\nextlist
                                        mark=o,mark=square,mark=triangle}}

\usepackage[maxbibnames=99, backend=bibtex]{biblatex}
\newcommand{\R}{\mathbb{R}}
\newcommand{\N}{\mathbb{N}}
\newcommand{\I}{\mathbb{I}}
\newcommand{\AI}{\operatorname{AI}}

\renewcommand{\vec}[1]{\text{vec}(#1)}
\newcommand{\GFlopps}{GFlop/s\xspace}

\newtheorem{remark}{Remark}
\newtheorem{example}{Example}
\title{Vectorized Parallel in Time methods for low-order discretizations with application to Porous Media problems}

\author[1]{Christian Engwer\thanks{christian.engwer@uni-muenster.de}}
\author[1]{Alexander Schell\thanks{Corresponding author: alexander.schell@uni-muenster.de}}
\author[1,2]{Nils-Arne Dreier\thanks{dreier@dkrz.de}}

\affil[1]{Applied Mathematics: Institute for Analysis and Numerics, University of Münster, Einsteinstr. 62, 48149 Münster, Germany}
\affil[2]{now at Deutsches Klimarechenzentrum, Bundesstraße 45a, 20146 Hamburg, Germany}

\date{}

\begin{document}

	\maketitle

	\begin{abstract}
		\noindent High order methods have shown great potential to overcome performance issues of simulations of partial differential equations (PDEs) on modern hardware, still many users stick to low-order, matrix-based simulations, in particular in porous media applications. Heterogeneous coefficients and low regularity of the solution are reasons not to employ high order discretizations.
		
		We present a new approach for the simulation of instationary PDEs that allows to partially mitigate the performance problems. By reformulating the original problem we derive a parallel in time time integrator that increases the arithmetic intensity and introduces additional structure into the problem. By this it helps accelerate matrix-based simulations on modern hardware architectures.
		
		Based on a system for multiple time steps we will formulate a matrix equation that can be solved using vectorized solvers like Block Krylov methods. The structure of this approach makes it applicable for a wide range of linear and nonlinear problems. In our numerical experiments we present some first results for three different PDEs, a \textit{linear convection-diffusion equation}, a \textit{nonlinear diffusion-reaction equation} and a realistic example based on the \textit{Richards' equation}.
	\end{abstract}


	\section{Introduction}\label{sec1}
In the context of porous media problems, simulations are often on a scale that requires significant computation power. Making good use of the available compute power on modern hardware architectures is thus of high importance.

As of today in porous media simulations low-order methods, like finite-volume or first-order mixed-FE methods, are still the dominant approach, e.g. \parencite{diersch2013feflow,koch2021dumux,rasmussen2021open,kolditz2012opengeosys,jung2017tough3}, and are usually following the classical way of a nonlinear Newton-type solver with an embedded iterative linear solver, where the latter is most often matrix based. One reason for the matrix based workflow is that it allows focussing on the modelling, while the solver \emph{simply works}. A second reason is the use of strong preconditioners, which often require a matrix to set up, e.g. ILU or AMG preconditioners.

As a consequence these methods are memory bandwidth limited and show a low arithmetic intensity~($\AI$), i.e. the ratio between floating point operations per transferred byte. This reduces their efficiency on modern hardware significantly.
  
The aim of this work is to offer a relatively lightweight approach to reformulate the existing methods in a way that improves hardware efficiency without modifying the spatial discretization. The central idea is to modify the time-integration method.

\paragraph*{Efficiency challenges}
Regarding efficient computations on modern hardware one thus faces two major problems. As the memory bandwidth is growing slower than the FLOP rate, low-order methods can not fully benefit from the hardware evolution. The second problem is an increasing concurrency. While many codes support the parallel simulation using multiple threads or processes, the data level parallelism, offered by SIMD ({Single Instruction Multiple Data}) and wide-SIMD units, requires sufficient structure in the numerical method. In the end, these problems hold true to CPUs just as well as GPUs.

The main computational costs are in the linear solver, which is usually bandwidth limited, due to the inherently low arithmetic intensity of the sparse matrix-vector product (\texttt{spMV}). For a double precision calculation two FLOP (one multiply and one add) are performed for two transferred \texttt{double} values (one matrix entry, one vector entry); the result entry is written multiple time, so we might neglect it here. All in all we obtain $\AI(\text{\texttt{spMV}})= \frac 1 8$, which means that the actual compute power is not used at all and on many modern computers we obtain only fraction of the peak-performance. The advertised peak  performance of an \textit{Intel Xeon Gold 6148 @ 2.40 GHz} is 1024~\GFlopps\footnote{\href{https://www.intel.com/content/www/us/en/content-details/840270/app-metrics-for-intel-microprocessors-intel-xeon-processor.html}{APP Metrics for Intel® Microprocessors - Intel® Xeon® Processor}}. Since it has 20 cores, the peak performance of a single core is 51.2~\GFlopps. Without AVX512 and FMA this reduces to 3.2~\GFlopps. Therefore, classical non-parallel code is bound by 3.2~\GFlopps and you loose $>$ 99\% of performance.

\paragraph*{Increasing the Arithmetic Intensity}
Not modifying the spatial discretization, two approaches to increase the  Arithmetic Intensity are possible.

One way to partially mitigate this problem is to use mixed-precision calculations, e.g. by applying the preconditioner in single precision and only calculate the residual in double precision. By this one is able to reduce the transfer costs, but still this increases $\AI$, and by this the actual performance, by less than a factor $2$. At the same time one has to be very careful to ensure that the reduced precision does neither impede the overall convergence and does not introduce stability issues. For an overview of mixed-precision methods and their analysis we refer to the review article by \cite{higham2022mixed}. In \parencite{grant2022perturbed} it is discussed how to apply mixed-precision algorithms to Runge-Kutta schemes.

An alternative approach is possible whenever a set of different linear systems with the same left-hand-side operator, but different right-hand-sides, needs to be solved. Examples are certain multi-scale methods, like MsFEM\,\parencite{hou1997multiscale} and LOD\,\parencite{engwer2019efficient}, certain inverse problems\,\parencite{lynn1967proposed,rullmann2009eeg} or the randomized construction of low-rank approximations\,\parencite{martinsson2020randomized}. Such a system can be reformulated to solve a single matrix equation instead of many linear systems with different right-hand-sides. The arising system takes the form 
\begin{equation}
    AX=B,
    \label{eq:basic-matrix-equation}
\end{equation} 
where $A \in \R^{N\times N}$ and $X,B\in \R^{N\times k}$ for some $k \in \N$, and can be solved using Block Krylov methods, which were first presented by \cite{o1980block}. They show better convergence behaviour, due to information  exchange between the columns. Recently Block Krylov methods got rediscovered to mitigate the communication overhead on high performance systems \parencite{grigori2016enlarged}. Their important feature is that they increase the arithmetic intensity and directly expose a vectorized structure \parencite{dreier2020strategies,bastian2020exa,lund2018new}. This is what the approach in this paper builds upon.

\paragraph*{Main Contribution and Structure of the paper}
The main contribution of this paper is to provide a new framework for the implicit solvers for instationary PDEs of the form
\begin{equation}
  \begin{aligned}
    M \dot{y}(t) + K y(t) &= b(t) \quad \text{for } t > 0, \\
    y(0) &= y_0,
  \end{aligned}
  \label{eq:problem}
\end{equation}
where $K \in \R^{N \times N}$ is the stiffness matrix, $b(t) \in \R^N$ a right-hand-side vector and mass matrix $M \in \R^{N \times N}$. In particular, such systems occur at the discretization of parabolic PDEs using the method of lines. 

The new approach picks up ideas from parallel in time methods (PinT). We transform the original problem into a space-time problem, considering multiple time steps at once. A similar approach was already presented by \cite{dekker2009partitioned}. In contrast to the traditional PinT paradigm we do not consider large time slaps and aim for a global parallelization over the time steps, but instead parallelize over consecutive time steps or Runge-Kutta stages to increase the node-level performance and employ vectorization. Vectorization is omnipresent on modern architectures via SIMD and wide-SIMD units, e.g. AVX and AVX-512 on x64 CPUs. To fully employ SIMD Units, additional structure is needed. We present a \emph{matrix splitting approach} to achieve this structure and formulate a matrix equation for solving the consecutive time steps. This matrix equation can then be solved by a vectorized Block Krylov method\,\parencite{dreier2020strategies}.

The method is fairly general and is especially well suited for stiffly accurate DIRK methods (\textit{diagonally implicit Runge-Kutta}) methods.

This paper is structured as follows. In Section~\ref{sec:VectorizedTimeIntegrators} we deduce the framework in an abstract way. In Section~\ref{sec:PerfomanceModelling} we explain the expected convergence behaviour with a simple performance model. Three numerical examples based on three different PDEs are presented in Section~\ref{sec:NumericalExperiments}. We conclude the work in Section~\ref{sec:Conclusion} and give a short outlook for further investigations.

	\section{Vectorized time integrators}\label{sec:VectorizedTimeIntegrators}

For the time discretization of problem~\eqref{eq:problem} 
consider the generic class of $m$~-~stage diagonally 
implicit Runge-Kutta methods (DIRK methods). 
These methods can be described by their Butcher tableau of the form
\begin{equation*}
    \begin{array}{c|c c c c c }
        c_1 & a_{1,1} & & & & \\
        c_2 & a_{2,1} & a_{2,2} & & & \\
        c_3 & a_{3,1} & a_{3,2} & a_{3,3} & & \\
        \vdots & \vdots & \vdots & \vdots & \ddots & \\
        c_m & a_{m,1} & a_{m,2} & a_{m,3} & \ldots & a_{m,m} \\
        \hline 
        & b_1 & b_2 & b_3 & \ldots & b_m       
    \end{array}
\end{equation*}
respectively the matrix $A\in \R^{m\times m}$ and the vectors $b,c\in\R^m$. 
We focus on \textit{stiffly accurate} methods that is methods 
satisfying the following property
\begin{equation*}
    b_i = a_{m,i} \quad \forall i = 1, \ldots , m.
\end{equation*}

Let $y^{(n,k)}\approx~y(t^{n}+c_k\tau)$ be the $k$-th stage of the 
$n$-th time step of the method and let $t^{(n,k)}~=~t^n+c_k\tau$ denote 
its corresponding time point. Then a stiffly accurate DIRK method to solve~\eqref{eq:problem} 
reads as
\begin{align}
    \begin{split}
        y^{(n,1)} &= y^{n} \\
        \tfrac{1}{\tau} M (y^{(n,k)} - y^{n}) + \sum_{j=1}^{k} a_{ij} K y^{(n,j)}
        &= \sum_{j=1}^{k} a_{i,j} b(t^{(n,j)}),
        \quad \forall k=1,\ldots,m \\
        y^{n+1} &= y^{(n,m)}
    \end{split}
    \label{eq:stiffly-accurate-DIRK}
\end{align}
with time step width $\tau>0$.

\subsection{Matrix equation formulation}

To formulate a matrix equation for all stages, one reformulates the system~\eqref{eq:stiffly-accurate-DIRK} into an all-at-once system, as typically done in \emph{PinT}-methods. Using the Kronecker product $\otimes$ the system for the $m$ stages of a stiffly accurate DIRK method can be written as
\begin{equation}
    \left[A_1 \otimes M + A_2 \otimes K\right] \vec{Y} = 
    \left[A_2 \otimes \I_N\right] \vec{B} + \vec{B_0},
    \label{eq:general-system-time-steps}
\end{equation}
where $Y,B \in \R^{N \times m}$ are defined as
\begin{equation*}
    \begin{array}{ccc}
        Y = (y^{(n,1)} , \ldots , y^{(n,m)}) & \mbox{and} & B = 
        (b^{(n,1)}, \ldots , b^{(n,m)})
    \end{array}
\end{equation*}
and $\vec{X} = (X_1^T,\ldots,X_s^T)^T \in \R^{(Nm)}$ denotes the stacked column vector of $X \in \R^{N \times m}$. 
The time stepping matrices $A_1,A_2 \in \R^{m \times m}$ and 
$B_0\in\R^{N\times m}$ depend on the particular DIRK method and are defined via 
the Butcher tableau
\begin{equation*}
    \begin{aligned}
        A_1 = \tfrac{1}{\tau}\I_m && , && A_2 = A && \mathrm{and} && 
        B_0 = \tfrac{1}{\tau} (M y^{n}, \ldots , M y^{n}).
    \end{aligned}
\end{equation*}

In the case of an explicit first stage, i.e. $a_{1,1}=0$, it is clear 
that $y^{(n,1)}~=~y^{n}$ and thus no linear system needs to be solved 
for this stage. Instead, $y^{(n,1)}$ can be moved to the right-hand-side directly.
Moreover, this effects the structure of the time stepping matrix $A_2$ and is well-suited for the \emph{matrix splitting approach} we present later.
We illustrate
this for general Runge-Kutta schemes. That is for schemes of the form
\begin{equation*}
    \begin{NiceArray}{c|c:ccc}
        0 & 0 & \Block{1-3}{} & & \\
        \hdottedline
        c_2 & a_{21}  & \Block{3-3}{A_\mathrm{impl}} & & \\
        \vdots & \vdots & & & \\
        c_m & a_{m1} & & & \\
        \hline
        & b_1 & b_2 & \cdots & b_m 
    \end{NiceArray}
\end{equation*}
with $A_\mathrm{impl} \in \R^{(m-1) \times (m-1)}$. Based on \eqref{eq:stiffly-accurate-DIRK} the system for the stages is given by
\begin{multline*}
    \left[ \tfrac{1}{\tau}
        \begin{pNiceArray}{c:ccc}
            1 & 0 & \cdots & 0 \\
            \hdottedline
            0 & \Block{3-3}{\I_{m-1}} & & \\
            \vdots & & & \\
            0 & & &
        \end{pNiceArray} \otimes M +
        \begin{pNiceArray}{c:ccc}
            1 & 0 & \cdots & 0 \\
            \hdottedline
            a_{21} & \Block{3-3}{A_\mathrm{impl}} & & \\
            \vdots & & & \\
            a_{m1} & & &
        \end{pNiceArray} \otimes K
    \right]
    \begin{pmatrix}
        y^{(n,1)} \\ y^{(n,2)} \\ \vdots \\ y^{(n,m)}
    \end{pmatrix} \\
    =
    \left[
        \begin{pNiceArray}{c:ccc}
            1 & 0 & \cdots & 0 \\
            \hdottedline
            a_{21} & \Block{3-3}{A_\mathrm{impl}} & & \\
            \vdots & & & \\
            a_{m1} & & &
        \end{pNiceArray} \otimes \I_N
    \right]
    \begin{pmatrix}
        b^{(n,1)} \\ b^{(n,2)} \\ \vdots \\ b^{(n,m)}
    \end{pmatrix} + \tfrac{1}{\tau}
    \begin{pmatrix}
        M y^n \\ M y^n \\ \vdots \\ M y^n
    \end{pmatrix}
\end{multline*}
Using the identity $y^{(n,1)} = y^{n}$ the corresponding first
row can be eliminated and the couplings to $y^{(n,1)}$ are moved
to the right-hand-side. This yields
\begin{multline*}
    \left[
        \tfrac{1}{\tau} \I_{m-1} \otimes M + A_\mathrm{impl} \otimes K
    \right]
    \begin{pmatrix}
        y^{(n,2)} \\ \vdots \\ y^{(n,m)}
    \end{pmatrix}
    \\
    = 
    \left[
        A_\mathrm{impl} \otimes \I_N
    \right]
    \begin{pmatrix}
        b^{(n,2)} \\ \vdots \\ b^{(n,m)}
    \end{pmatrix}
    + \tfrac{1}{\tau}
    \begin{pmatrix}
        M y^n \\ \vdots \\ M y^n
    \end{pmatrix}
    +
    \begin{pmatrix}
        a_{21} (b^n - K y^n) \\ \vdots \\ a_{m1} (b^n - K y^n)
    \end{pmatrix}
\end{multline*}
above~\eqref{eq:general-system-time-steps} for a DIRK method
which is exactly of the form described
of one stage less and couplings to the eliminated stage represented by additional contributions
in $B_0$.

\begin{example}
    Eliminating the first stage of the $\Theta$-method this yields its well known form
    \begin{equation*}
        \left(\tfrac{1}{\tau} M + \Theta K \right) y^{n+1} = \Theta b^{n+1} + \tfrac{1}{\tau} M y^n + (1 - \Theta) (b^n - K y^n),
    \end{equation*}
    which will be used in the later numerical experiments in Section~\ref{sec:NumericalExperiments}.
\end{example}

These concepts can easily be extended such that not only $m$ stages of
a single time step are considered, but also $s$ consecutive time steps, each consisting of 
$m$ stages. This is directly achieved by extending the matrices 
$A_1,A_2$ to block matrices in $\R^{(s\cdot m) \times (s \cdot m)}$ and $B_0$ 
to $\R^{N \times (s \cdot m)}$, leading to the following matrices
\begin{equation*}
    \begin{aligned}
        A_1 = \tfrac{1}{\tau}\begin{pmatrix}
            \I_m & & \\
            & \ddots & \\
            & & \I_m
        \end{pmatrix},
        && 
        A_2 = \begin{pmatrix}
            A & & \\
            & \ddots & \\
            & & A
        \end{pmatrix}
        && \mathrm{and} 
    \end{aligned}
\end{equation*}
\begin{equation*}
    B_0 = \tfrac{1}{\tau} 
    \begin{pmatrix}
        \underbrace{M y^{n} , \,\ldots\, , M y^{n}}_{m} , \; \ldots \; ,
        \underbrace{M y^{n+s-1} , \, \ldots \, , M y^{n+s-1}}_{m}
    \end{pmatrix},
\end{equation*}
with potentially additional terms in $B_0$ if an explicit stage was 
eliminated before. Note that the $y^{n+i} , i=1,\ldots,s-1$ are in general
unknown, but for the stiffly
accurate methods it holds that
\begin{equation*}
    y^{n+i} = y^{(n+i-1,m)} , \quad i=1,\ldots,s-1
\end{equation*}
and hence the corresponding dependencies are already included in the left-hand-side operator. The matrices therefore will look as follows 

\begin{equation}
    \begin{array}{lll}
        \multicolumn{3}{c}{
            A_1 =  
            \begin{pNiceMatrix}[columns-width=5mm]
                \Block[borders={bottom,right}]{3-3}{\I_m} & & & & & & & & \\
                & & & & & & & &\\
                & & & & & & & &\\
                & & -1 & \Block[borders={bottom,right,top,left}]{3-3}{\I_m} & & & & &\\
                & &  \vdots & & & & & &\\
                & & -1 & & & & & &\\    
                & & & & & -1 &\Block[borders={top,left}]{3-3}{\ddots} & &\\
                & & & & & \vdots & & &\\
                & & & & & -1 & & &\\
            \end{pNiceMatrix},
        } \\[2.5cm]
        A_2 = \begin{pNiceMatrix}[columns-width=5mm]
            \Block[borders={bottom,right}]{1-1}{A} & & \\
            & \Block[borders={left,top,right,bottom}]{1-1}{A} & \\
            & & \Block[borders={left,top}]{1-1}{\ddots}
        \end{pNiceMatrix}
        & \mathrm{and} &
        B_0 = \tfrac{1}{\tau} 
        \begin{pmatrix}
            \underbrace{M y^{n} , \,\ldots\, , M y^{n}}_{m} , 0 , \, \ldots
        \end{pmatrix}.
    \end{array}
    \label{eq:time-stepping-matrices-multiple-parallel-time-steps}
\end{equation}

If an explicit stage was eliminated, $A_2$ will as well carry the additional 
entries $a_{2,1},\ldots,a_{m,1}$ in the same positions as the $-1$ entries 
in $A_1$ and $B_0$ will also contain the additional terms in the first 
$m$ columns. In Example~\ref{ex:time-stepping-matrices} we show how these 
matrices $A_1,A_2$ and $B_0$ will look like for the
$\Theta$-method.

From now on, for simplicity we assume that only a single time step consisting 
of $m$ stages is considered, but all concepts can be canonically extended to 
the more complex case of multiple time steps. Later, the main ideas are 
shown by the use of some examples.  

\begin{example}{(Time stepping matrices for $\Theta$-method)}
    Considering $s$ consecutive time steps of the $\Theta$-method
    after eliminating the explicit stage the time stepping matrices
    $A_1,A_2 \in \R^{s \times s}$ and the right-hand-side $B_0 \in
    \R^{N \times s}$ are given by
    \label{ex:time-stepping-matrices}
    \begin{equation*}
        \begin{array}{ll}
            A_1= \tfrac{1}{\tau}
            \begin{pmatrix}
                1 & & &  \\
                -1 & 1 & &  \\
                & \ddots & \ddots & \\
                & & -1 & 1
            \end{pmatrix},\quad
            &
            A_2= 
            \begin{pmatrix}
                \Theta & & &  \\
                1-\Theta & \Theta & &  \\
                & \ddots & \ddots & \\
                & & 1-\Theta & \Theta
            \end{pmatrix},
            \\[1.5cm]
            \multicolumn{2}{l}{
            B_0 = 
            \begin{pmatrix}
                \tfrac{1}{\tau} M y^{n} -(1-\Theta) (Ky^{n} - b^{n}) , & 0, & \ldots, &0
            \end{pmatrix}.
            }
        \end{array}
    \end{equation*}
\end{example}
To apply vectorized Block Krylov methods to the resulting system it must be ensured
that~\eqref{eq:general-system-time-steps} can be transformed into a 
matrix equation.
Therefore, for matrices $A \in \R^{N \times N}$ and $X,B \in R^{N \times m}$, recall the general identity
\begin{align*}
    \begin{pmatrix}
        A & & \\
        & \ddots & \\
        & & A
    \end{pmatrix}
    \begin{pmatrix}
        X_1^T \\ \vdots \\ X_m^T
    \end{pmatrix}
    & =
    \begin{pmatrix}
        B_1^T \\ \vdots \\ B_m^T
    \end{pmatrix}\\
  \Leftrightarrow \left( \I_m \otimes A \right) \vec{X}
  &= \vec{B} \\
  \Leftrightarrow AX &= B.
\end{align*}
Thus, instead of solving multiple linear system with the same matrix $A$ we can solve a matrix equation of the form \eqref{eq:basic-matrix-equation}.
To obtain the desired structure for the left-hand-side operator 
we propose a \emph{matrix~splitting~approach} that allows us to decouple the
stages in the operator, moving any 
off-diagonal entry to 
the right-hand-side. Split the matrices  $A_1$ and $A_2$ into their 
diagonal and off-diagonal parts
\begin{equation}
    A_1 = \tfrac{1}{\tau} \I_m + \left( A_1 - \tfrac{1}{\tau} \I_m\right) \quad \text{and} 
    \quad A_2 = \beta \I_m + \left( A_2 - \beta \I_m\right),
    \label{eq:matrix-splitting}
\end{equation} 
with $\beta = A_{2_{(0,0)}}$ the entry in the 
top left corner of the matrix and move all off-diagonal parts to the 
right-hand-side of the equation. Equation~\eqref{eq:general-system-time-steps} 
then can be written as
\begin{multline*}
    \left[\tfrac{1}{\tau} \I_m \otimes M + \beta \I_m \otimes K \right] \vec{Y} \\
     = - \left[(A_1 - \tfrac{1}{\tau} \I_m) \otimes M + (A_2 - \beta \I_m) 
     \otimes K \right] \vec{Y} + \left[A_2 \otimes \I_N\right] \vec{B} + \vec{B_0}  ,
\end{multline*}
which has an equivalent matrix equation formulation
\begin{equation*}
    \left(\tfrac{1}{\tau} M + \beta K\right) Y = - M Y (A_1 - \tfrac{1}{\tau} \I_m)^T - 
    K Y (A_2 - \beta \I_m)^T + B A_2^T + B_0.
\end{equation*}

Due to the splitting, the derived matrix equation 
has a right-hand-side depending on the solution $Y$. Therefore, it cannot 
be solved directly with a vectorized Block Krylov solver but has to be 
solved iteratively, e.g. using a fix-point iteration. During this iteration 
the actual matrix equation then can be solved using a vectorized solver
\begin{align}
    \begin{split}
    \left(\tfrac{1}{\tau} M + \beta K\right) Y^j & = - M Y^{j-1} (A_1 - \tfrac{1}{\tau} \I_m)^T - 
    K Y^{j-1} (A_2 - \beta \I_m)^T + B A_2^T + B_0
    \end{split}
    \label{eq:matrix-fixed-point-equation}
\end{align}

\begin{remark}
    Let the mass matrix $M$ be a scaled identity and the left-hand-side 
    operator take the form
    \begin{equation*}
        D \otimes M + \beta \I_m \otimes K,
    \end{equation*}
    with $D \in \R^{m \times m}$ a diagonal matrix. In that case one can
    apply vectorized Block Krylov methods, even though the problem is
    not in the form of a matrix equation. The reason is that the Block Krylov 
    space of $(D \otimes M + \beta \I_m \otimes K)$ is the same as that 
    of $(D_{00} M + \beta K)$.
\end{remark}

\begin{remark}
    If one considers the matrix splitting~\eqref{eq:matrix-splitting} e.g. for the time stepping matrices of the $\Theta$-method from Example~\ref{ex:time-stepping-matrices}, one realizes that eliminating the explicit stage leads to a smaller splitting error since all zero diagonal entries were removed. 
\end{remark}

\subsection{Extension to the nonlinear case}\label{subsec:extension-to-nonlinear-case}

For nonlinear problems, the stiffness matrix has to be replaced by a nonlinear operator $\mathcal{K}$. The nonlinear version of the semi-discretized system~\eqref{eq:problem} then reads
\begin{equation*}
    M \dot{y}(t) + \mathcal{K}(y(t)) = b(t).
\end{equation*}
After linearizing $\mathcal{K}$ at $y(t)$ this corresponds to solve
\begin{equation*}
    M \dot{y}(t) + K[y(t)]\, y(t) = b(t),
\end{equation*}
with stiffness matrix $K[y(t)] \in R^{N \times N}$ or equivalently with $f(t,y(t)):=K[y(t)]\,y(t)-b(t)$
\begin{equation*}
    M \dot{y}(t) + f(t,y(t)) = 0.
\end{equation*}
Thus, the stiffness matrix will differ between the time steps since different linearization points are used. To formulate the matrix equation in this case we
introduce another splitting, which ensures
that the left-hand-side operator is the same for all time
steps.
Similar to \cite{dekker2009partitioned} we use an approximate
$\tilde{K}$ for every time steps. Whilst they propose to define
$\tilde{K}$ as the average over all stiffness matrices, we will
use the stiffness matrix of the first time step.
In the same fashion as before one can 
formulate a system for the $m$ parallel stages of a stiffly accurate 
DIRK method using the Kronecker product
\begin{equation}
    \left[A_1 \otimes M + (A_2 \otimes \I_N) K \right] \text{vec}(Y) 
    = \left[ A_2 \otimes \I_N\right] \vec{B} + \vec{B_0},
    \label{eq:system-for-time-steps-nonlinear-case}
\end{equation}
with $K$ defined as
\begin{equation*}
    K := 
    \begin{pmatrix}
        K[y^{(n,1)}] & & \\
        & \ddots & \\
        & & K[y^{(n,m)}]
    \end{pmatrix}
    :=
    \begin{pmatrix}
        K^{(n,1)} & & \\
        & \ddots & \\
        & & K^{(n,m)}
    \end{pmatrix}.
\end{equation*}

Since equation~\eqref{eq:system-for-time-steps-nonlinear-case} needs to be transformed
into a matrix equation as done in the linear case, as before a splitting has to be 
applied. Additionally, one needs to ensure that the same 
stiffness matrix, i.e. the same linearization is used for every stage. 
Therefore, similar to \cite{dekker2009partitioned} we use an approximate 
$\tilde{K}$ for all stages. While they suggest constructing $\tilde{K}$ 
by averaging, we propose to use the Jacobian of the first stage, similar 
to the reuse of the Jacobian in simplified Newton methods.
For each stage or time step we split the individual
Jacobian matrices $K^{(n,i)}$ as follows
\begin{equation*}
    K^{(n,i)} = K^{(n,1)} + \tilde{K}^{(n,i)}.
\end{equation*}

With the abbreviated notation
$
    \tilde{K} := \text{\footnotesize$
    \begin{pmatrix}
        \tilde{K}^{(n,1)} \hspace*{-1ex}& & \\[-1ex]
        & \ddots & \\[-1ex]
        & & \tilde{K}^{(n,m)}
    \end{pmatrix}$}
$
one gets
the following identity
$
K = (\I_{m} \otimes K^{(n,1)}) + \tilde{K},
$
which allows for a split reformulation of
system~\eqref{eq:system-for-time-steps-nonlinear-case} as follows
\begin{align}
    \begin{split}
        \Big[A_1 \otimes M + \underbrace{(A_2 \otimes \I_N)
        (\I_{m} \otimes K^{(n,1)})}_{=(A_2 \otimes K^{(n,1)})} \Big] 
        \text{vec}(Y) = & - \Big[(A_2 \otimes \I_N) \tilde{K}\Big] \vec{Y} 
        \\[-1ex]
        & + \Big[ A_2 \otimes \I_N\Big] \vec{B} + \vec{B_0}.
    \end{split}
\end{align}

Up to the additional correction term of the Jacobian on the right-hand-side,
this is the same formulation as in the linear 
case~\eqref{eq:general-system-time-steps}. Hence, from now on one can 
apply the same splitting as done in the linear case to obtain the matrix 
equation formulation
\begin{align}
    \begin{split}
        \left(\tfrac{1}{\tau} M + \beta K^{(n,1)}\right) Y = & - \left(\tilde{K}^{(n,1)} 
        y^{(n,1)},\ldots , \tilde{K}^{(n,m)} y^{(n,m)}\right) A_2^T \\
        & - M Y (A_1 - \tfrac{1}{\tau} \I_{m})^T \\
        & - K^{(n,1)} Y (A_2 - \beta \I_{m})^T + B A_2^T + B_0.
    \end{split}
    \label{eq:matrix-equation-nonlinear-case}
\end{align}

As in the linear case, equation~\eqref{eq:matrix-equation-nonlinear-case} has 
a solution dependent right-hand-side and thus needs to be solved iteratively,
but instead of solving the matrix equation exactly, we will fuse the fix-point 
iteration with the outer Newton update. We use one iteration of the 
fix-point problem as the search direction in the outer Newton solver. 
Due to the mostly explicit structure of the DIRK methods, i.e. stages only couple 
with previous stages, there is no splitting needed for the left-hand-side 
operator of the first stage when choosing the Jacobian of the first stage as the approximate Jacobian 
for all stages. Thus, it is to be expected that the first stage converges fastest, and 
we employ a so-called \emph{pipelining strategy}, shifting to the next 
stage, whenever the first stage has converged. As a consequence, after shifting, 
the first row in system~\eqref{eq:system-for-time-steps-nonlinear-case} 
will then correspond to the next stage and will not necessarily belong to 
the first stage of a time step. Hence, the matrices $A_1,A_2,B$ and $B_0$ 
actually depend on the stage considered for the top row. We indicate this 
by an additional index $A_1^{(n,k)},A_2^{(n,k)},B^{(n,k)}$ and $B^{(n,k)}_0$. 
To simplify the consecutive numeration of the stages in that case define 
the index mapping 
\begin{align*}
    \begin{split}
    \iota : & \N \times \N \to \N \times \{ 1 ,\ldots, m\} \\
    & (n,k) \mapsto \left( n + \left\lfloor \frac{k-1}{m} \right\rfloor , k~\bmod m \right)
    \end{split}
\end{align*}

Using this mapping denote the $m$ successive stages starting from $y^{(n,k)}$ 
by $y^{\iota(n,k)},\ldots,y^{\iota(n,k+m-1)}$. Especially, this means 
that some considered stages belong to different time steps, hence 
we are in a similar situation as described in~\eqref{eq:time-stepping-matrices-multiple-parallel-time-steps}. 
Given the stages $y^{(n,1)},\ldots,y^{(n,k-1)}$ the matrix 
equation~\eqref{eq:matrix-equation-nonlinear-case} for 
$Y^{(n,k)}~=~\left(y^{\iota(n,k)},\ldots,y^{\iota(n,k+m-1)}\right)$ formulates as
\begin{align}
    \begin{split}
        \Big(\tfrac{1}{\tau} M + \beta^{(n,k)}& K^{\iota(n,k)}\Big) Y^{(n,k)} \\
        = & - (\tilde{K}^{\iota(n,k)} 
        y^{\iota(n,k)},\ldots , \tilde{K}^{\iota(n,k+m-1)} y^{\iota(n,k+m-1)}) A_2^{(n,k)^T} \\
        & - M Y^{(n,k)} (A_1^{(n,k)} - \tfrac{1}{\tau} \I_{m})^T - K^{\iota(n,k)} Y^{(n,k)} (A_2^{(n,k)} - \beta^{(n,k)} 
        \I_{m})^T \\
        & + B^{(n,k)} A_2^{(n,k)^T} + B_0^{(n,k)},
    \end{split}
    \label{eq:matrix-equation-pipelining}
\end{align}
where $A_1^{(n,k)}$ and $A_2^{(n,k)} \in \R^{m \times m}$ the $m\times m$ 
submatrices from~\eqref{eq:time-stepping-matrices-multiple-parallel-time-steps} 
starting at the index $(k,k)$. The $i$~-th column of the right-hand-side 
$B_0^{(n,k)}$ is given by 
\begin{equation*}
    B_{0_i}^{(n,k)} = \begin{cases*}
        \sum_{j=1}^{k-1} a_{i,j} f(t^n + \tau c_j , y^{(n,j)}) + 
        \tfrac{1}{\tau}M y^{n} & if $i+k \leq m+1$ \\
        0 & else.
    \end{cases*}
\end{equation*}

Taking a closer look at the matrix equation~\eqref{eq:matrix-equation-pipelining},
one observes that the residual of this equation is given by
\begin{equation*}
    \mathcal{R}(Y^{(n,k)}) := M Y^{(n,k)} A_1^{(n,k)^T} - \mathcal{F}(t1^T + \tau c , Y^{(n,k)}) A_2^{(n,k)^T} - B_0^{(n,k)},
\end{equation*}
where $\mathcal{F}: \R^{m} \times \R^{N \times m} \to \R^{N \times m}$ 
is the vectorization of
$f(t,y(t)):= K[y(t)] y(t) - b(t)$ from~\eqref{eq:problem}, given by
\begin{equation*}
    \mathcal{F}(T,Y) = \left(f(T_1,Y_1),\ldots,f(T_m,Y_{m})\right) \in \R^{N \times m}.
\end{equation*}

The pseudocode is presented in Algorithm~\ref{alg:Block Krylov-Time-Stepping-Pipelining}
and is described in the following. The algorithm starts by initializing $Y^{(0,0)}$. An outer loop loops over all 
time steps and an inner loop loops over all stages, i.e. handles the 
pipelining and shifts the stages to the left using the function 
\lstinline{shift} in line~\ref{algline:shift}. 

\begin{equation*}
    \operatorname{shift}\left((Y_1, \ldots, Y_{s-1} , Y_s)\right) = 
    (Y_2 , \ldots , Y_s, Y_s )
    \label{eq:shift}
\end{equation*}

Next, the Newton iteration is started, where the residual and the Jacobian 
are computed as above and finally in line~\ref{algline:solve} the system is solved using a vectorized 
solver. 

\begin{algorithm}
    \caption{Block Krylov Time Stepping with Pipelining}
    \label{alg:Block Krylov-Time-Stepping-Pipelining}
    \begin{algorithmic}[1]
        \State $Y^{(0,0)} = (y_0,\ldots,y_0)$
        \For{time steps $n = 0,\ldots,$}
            \For{stages $k=1,\ldots,m$}
                \State $Y^{(n,k),0} = \operatorname{shift}(Y^{(n,k-1)})$ \label{algline:shift}
                \For{$i=0,\ldots$}
                    \State $R \gets \mathcal{R}(Y^{(n,k),i})$
                    \If{$\| R_0 \| < \varepsilon$}
                        \State $Y^{(n,k)} = Y^{(n,k),i}$ 
                        \State break                    
                    \EndIf
                    \State $J \gets \left(\tfrac{1}{\tau}M + \alpha_{k,k} Df(t^n + \tau c_k , Y^{(n,k),i}_0)\right)$
                    \State solve $JV=R$ \label{algline:solve}
                    \State $Y^{(n,k),i+1} = Y^{(n,k),i} - V$
                \EndFor
            \EndFor
        \EndFor
    \end{algorithmic}
\end{algorithm}

\begin{remark}[Fully-implicit Runge-Kutta methods]
    All reformulations shown above that lead to the matrix equation can also be done
    for fully-implicit Runge-Kutta schemes. The main difference is that for those
    methods there will be a splitting error in the first stage, hence one cannot
    guarantee the first stage to converge the fastest which contradicts the
    assumptions made for the \emph{pipelining}. In that case one should
    use an average Jacobian as proposed by \cite{dekker2009partitioned}.
\end{remark}

\begin{remark}[Additional Nonlinearity]
    In the case that the underlying problem has an additional nonlinearity 
    in the time derivative, which is the case e.g. for Richards' equation (see Section~\ref{sec:NumericalExperiments}),
    not only the stiffness matrix but also the mass matrix will be
    time dependent. As a consequence it is necessary to apply the same
    splitting concepts also to the mass matrix, resulting in an
    additional splitting error coming from the inexact linearization.
\end{remark}

\begin{remark}[Linearization]
    In some cases it might be useful to use a different linearization 
    than the classical Newton method for the spatial part (or the temporal part).
\end{remark}

	\section{Performance analysis}\label{sec:PerfomanceModelling}
A widely used tool to study the performance of an algorithm is 
the \emph{arithmetic intensity} ($\AI$), which measures the ratio of
computational operations to memory traffic and is defined as the number
of \emph{floating-point operations} (FLOPs) per data transferred between
memory and processor (measured in \emph{bytes}).
A main goal of our approach is to make use of vectorized Krylov
solvers to increase the performance by increasing the arithmetic
intensity. In this section we present a short performance analysis for
a vectorized Block Krylov method solving
\begin{equation*}
    A X = B, \quad \mbox{with } B,X \in \R^{N \times k}.
\end{equation*}

A main ingredient of an iteration of such a method is the \textsc{BOP}
kernel, i.e. applying the operator on a block vector
\begin{equation*}
    Y \mapsfrom A X.
\end{equation*}
Compared to the standard sparse matrix-vector product $(\AI(\texttt{spMV})=\frac{1}{8})$
this has an increased arithmetic intensity since the matrix only has to 
be loaded once if we make use of SIMD.

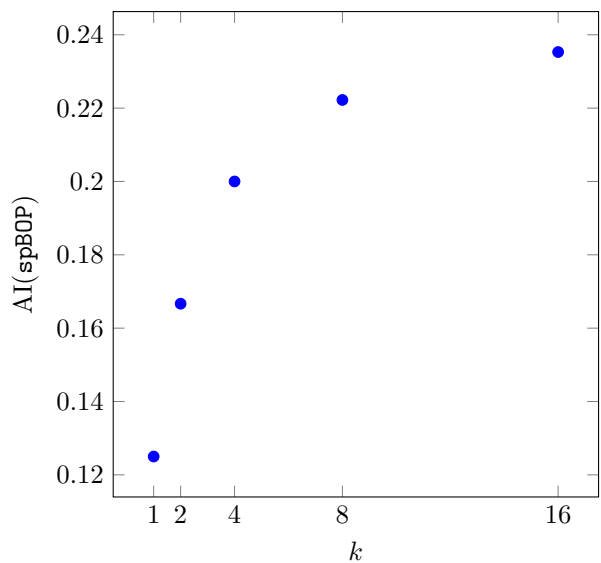
\begin{figure}[h]
    \centering
    \begin{tikzpicture}
        \begin{axis}
            [
            samples at={1,2,4,8,16},
            xtick = {1,2,4,8,16},
            log basis x = 2,
            xlabel = $k$,
            ylabel = $\AI(\texttt{spBOP})$,
            width=0.5\textwidth,
            height=0.5\textwidth,
            ]
            \addplot[only marks,mark=*,color=blue] {x/(4*(x+1))};
        \end{axis}
    \end{tikzpicture}
    \caption{Arithmetic intensity of the sparse \textsc{BOP} kernel for different 
    block sizes.}
    \label{fig:AIBOP}
\end{figure}

Consider a 
block-vector $X\in\R^{N\times k}$ with $k$ columns. For a double precision 
calculation two $2\cdot k$ FLOP (one multiply and one add for each 
column) are performed for $k+1$ transferred \textsc{double} values (one 
matrix entry and $k$ block-vector entries). As before, neglect the 
result entry here. Thus, we obtain an arithmetic intensity of the sparse 
\textsc{BOP} kernel $\AI(\texttt{spBOP}) = \frac{k}{4(k+1)}$.
Figure~\ref{fig:AIBOP} shows the arithmetic intensity of the sparse 
\textsc{BOP} kernel for different values of $k$. One sees that for $k=16$ the arithmetic intensity increased by a factor of about 1.8.

To analyze the performance of the \textsc{BOP} we consider a simplified
variant of the ECM model presented by \cite{hofmann2019bridging} like
in \cite{dreier2020strategies}, where they assume that the runtime of
a kernel is bounded by three factors:

\begin{itemize}
    \item $T_\text{comp} = \frac{\omega}{\text{peakflops}}$: The time 
    the processor needs to perform the necessary number of floating point 
    operations, with $\omega$ denoting the number of floating point 
    operations of the kernel.
    \item $T_\text{mem} = \frac{\beta}{\text{memory bandwidth}}$: The 
    time to transfer the data from the main memory to the L1-cache, with 
    the amount of data that needs to be transferred in the kernel denoted 
    by $\beta$.
    \item $T_\text{reg} = \frac{\tau}{\text{register bandwidth}}$: The 
    time to transfer data between L1-cache and registers. Here $\tau$ 
    is the amount of data that needs to be transferred in the kernel.
\end{itemize}
The effective execution time is then given by
\begin{equation*}
    T = \max{(T_\text{comp} , T_\text{mem} , T_\text{reg})}
\end{equation*}
We denote by $\omega$ the number of floating points operations for
this kernel, by $\beta$ the amount of data loaded from main memory and
the number of data transfers between registers and L1-cache is given
by $\tau$. The number of nonzero entries of $A$ is denoted by $z$. In
the case of the \textsc{BOP} kernel these characteristics are given by
\begin{align*}
    \omega & = 2 k z\\
    \beta & = 2z + 2kz\\
    \tau & = z(2+2k)
\end{align*}

\newcommand{\flops}[3]{2*(#1*#2*#3)}
\newcommand{\LD}[3]{8*#3*#2*(8+2*#1)} 
\newcommand{\ST}[3]{8*(#1*#3)}
\newcommand{\mem}[3]{8*#3*(2*#2+2*#1)}

\begin{figure}[h]
    \centering
    \begin{subfigure}[t]{0.45\textwidth}
        \centering
        \includegraphics[width=\linewidth,height=\linewidth]{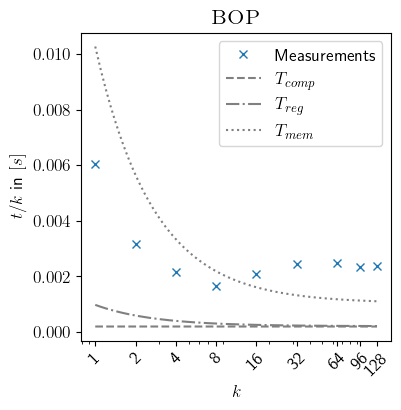}    
        \caption{Measured and expected runtime per right-hand-side of 
        the \textsc{BOP} kernel.}
        \label{fig:bopbenchmark}
    \end{subfigure}
    ~
    \begin{subfigure}[t]{0.45\textwidth}
        \centering
        \includegraphics[width=\linewidth,height=\linewidth]{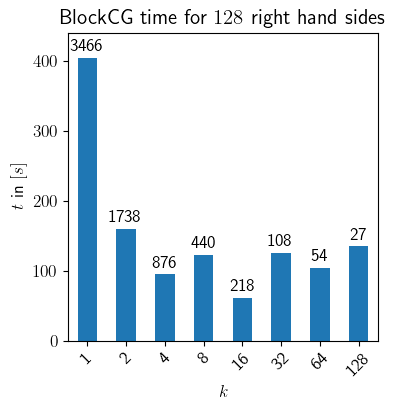}
        \caption{Time needed by the Block-CG to solve for $128$ 
        right-hand-sides considering $k$ right-hand-side in parallel. Numbers 
        on top of the bars denote the number if iterations.}
        \label{fig:bcgbenchmark}
    \end{subfigure}
    \caption{Microbenchmark for the \textsc{BOP} kernel (left) and the 
    Block-CG method (right).}
    \label{fig:MicrobenchmarkBOPKernel}
\end{figure}

Figure~\ref{fig:bopbenchmark} shows the measured runtimes
compared to the expected ones. All computations are done on an
\emph{Intel Xeon Gold 6148} on a single core. Its theoretical
peak performance is
38.4 \GFlopps, the memory bandwidth is 6.39 Gbyte/s and
the L1-cache bandwidth is
286.1 Gbyte/s. The microbenchmark for the 
\textsc{BOP} kernel shows that as expected the runtime per right-hand-side decreases 
rapidly for small $k$. For larger $k$ the runtime per right-hand-side 
increases slightly. A reason for this might be that for larger $k$ one 
row of a block vector occupies more space in the caches. Thus, fewer rows 
can be cached resulting in a longer runtime per right-hand-side.
Figure~\ref{fig:bcgbenchmark} further shows the time needed to
solve an elliptic PDE for $128$
right-hand-sides using a Block-CG Krylov solver with an AMG preconditioner.
The exact problem considered is the stiffness matrix of the linearized
Richards' equation (see Section~\ref{subsec:Richards}).
The required accuracy is a defect reduction of $10^{-8}$ for all right-hand-sides.
One clearly observes a significant reduction in computation time. In
this particular case the sweet spot is for $k=16$. Even though fewer 
iterations are needed when increasing the block size, we do not get faster 
for $k$ larger than $16$, since each row of a block vector occupies more 
space in the cache like we observed for the \textsc{BOP} kernel.

\begin{remark}[Nonlinearity]
  \label{remark:non-linear-convergence}
  In addition to discussed benefits in the linear case,
  we expect a secondary effect in the nonlinear case.

  Solving \eqref{eq:matrix-equation-pipelining} with Newton's method, the exact Jacobian 
  is used for the first stage and an approximate Jacobian for any
  further stages. This means, effectively we perform an exact Newton iteration for the first stage and 
  an inexact Newton iteration for the rest,
  which introduces a linearization error. Still the inherent splitting error
  usually dominates the convergence of the outer iteration, i.e. the
  nonlinear solver. While we don't actually solve the exact
  nonlinear problem on later stages, still we compute an improved
  initial guess.
  This means that after shifting, when a stage become the first
  one, the initial guess is better than in the sequential case.
  Often this effect becomes even more prominent for later time steps, resulting in less Newton iterations per stage and an 
  additionally reduced number of linear solver iterations. This can lead 
  to a behaviour where after some initial iterations every stage 
  converges after a single Newton iteration after it has been shifted to 
  the front by the pipelining.

This clearly depends on the dynamic of the problem we look at. Since 
we use approximate Jacobians, this behaviour depends on the dynamics of 
the nonlinear operator $f$, i.e. its derivative
\begin{equation*}
    \frac{d}{dt} f(t , y(t)).
\end{equation*}
The more the nonlinear operator changes over time, the worse the 
approximation error gets we make by choosing the wrong linearization point 
for the Jacobian.
\end{remark}

	\section{Numerical Experiments}\label{sec:NumericalExperiments}
\begin{minipage}{1.0\linewidth} 
  In the numerical experiments we test the method on a series of
  different problems with growing complexity:
  \begin{enumerate}
  \item Linear Convection-Diffusion
  \item Nonlinear Diffusion-Reaction
  \item Richards' equation
  \end{enumerate}
  \smallskip
\end{minipage}

We implemented a prototype of the pipelining approach presented in
section~\ref{sec:VectorizedTimeIntegrators} in \textsc{DUNE}
\parencite{bastian2008generic,bastian2021dune,dune:pdelab} where we made use of \textsc{DUNE}'s SIMD abstraction to
benefit from the SIMD capabilities of the hardware.

\begin{figure}
  \centering
  \includegraphics[width=0.8\linewidth]{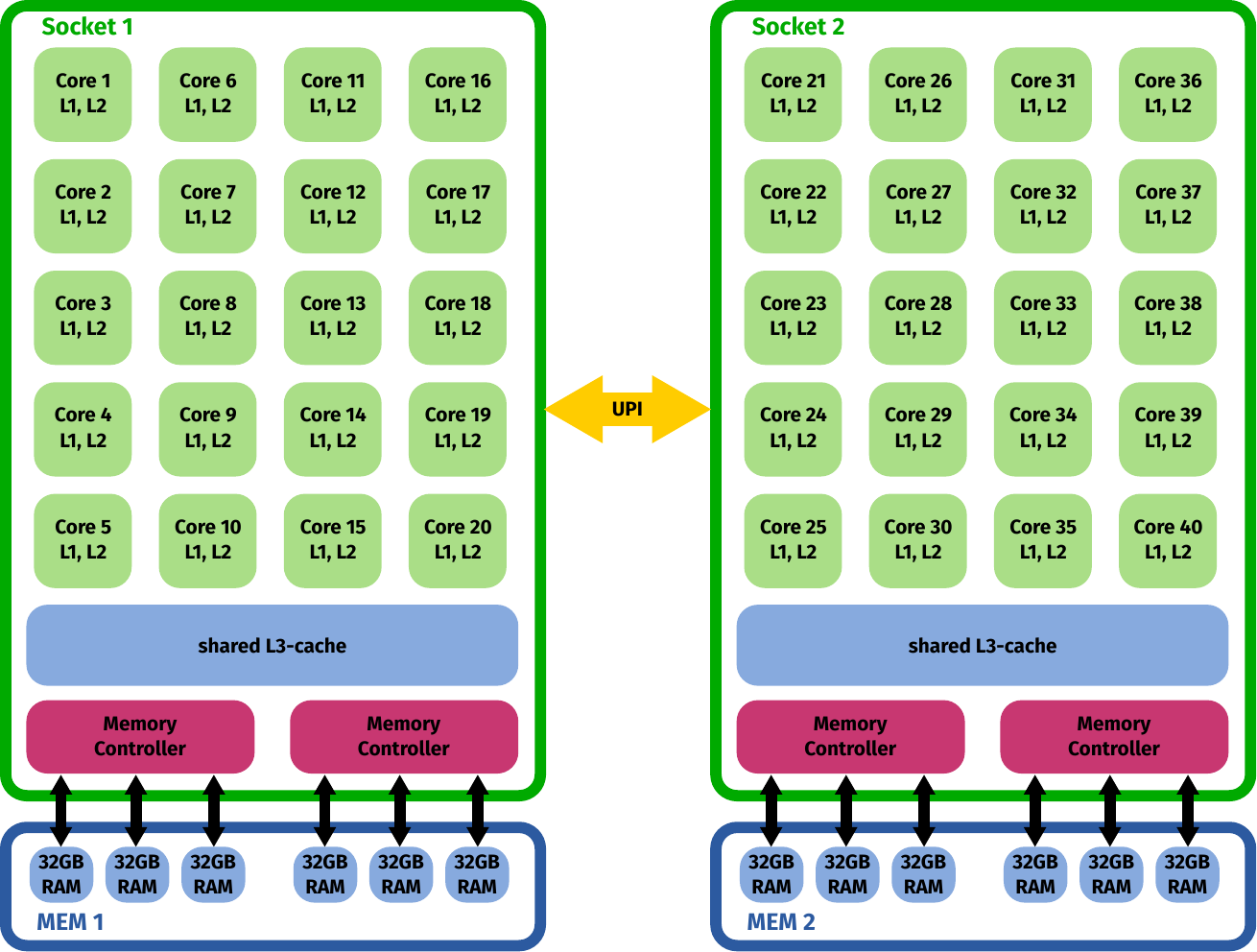}
  \caption{Sketch of the hardware topology of the \emph{Intel(R) Xeon(R)
        Gold 6148 CPU @ 2.40GHz} used for the numerical
      experiments. 2 sockets with 20 cores each yields a NUMA
      layout. Each core has separate L1- and L2-cache, whereas the
      cores in each socket share one 28 MB L3-cache.}
    \label{fig:HardwareTopology}
\end{figure}

All calculations were performed on an \emph{Intel Xeon Gold 6148 CPU
  @ 2.40GHz} (Intel Skylake) with
\texttt{gcc} (11.4.0) and
\emph{AVX512} support enabled (\texttt{-march=skylake-avx512}).
The exact hardware topology of the machine
is depicted in Figure~\ref{fig:HardwareTopology}. It consists of two sockets
with 20 cores each that have a separate L1- and L2-cache but share the same L3-cache.
We ran all of our tests on a single socket, executing 20 identical
processes via MPI in order to occupy all cores
\begin{center}
    \begin{lstlisting}{language=bash}
        mpirun -n 20 --bind-to core --map-by core <executable file>        
    \end{lstlisting}
\end{center}
This mimics the situation in a productive parallel run, where all
cores share a common L3-cache, but compete for the same 6 memory interfaces.

\subsection{Linear Convection-Diffusion}

The first example considered is a linear convection-diffusion equation
\begin{align*}
    \partial_t u - \nabla\cdot(-A(x) \nabla u + b(x) u) + c(x) u &=0 & \mbox{ in } \Omega \times (0,T],  \\
    u(t,x) &= g(t,x) & \mbox{on } \partial\Omega_D \times (0,T],  \\
    (b(x)u-(A\nabla u)) \cdot n &= 0 & \mbox{ on } \partial\Omega_O \times (0,T], \\
    u(0,x) & = u_0(x) & \mbox{in } \Omega,
\end{align*}
where the computational domain is given by $\Omega = [0,1]^2$, 
the Dirichlet boundary is given by
$\partial\Omega_D =\Gamma_1 \cup \Gamma_2 = \{0\} \times [0,1] \cup
[0,1] \times \{1\}$ and outflow boundary
$\partial\Omega_O = \partial\Omega \setminus \partial\Omega_D$. We
choose the permeability tensor $A(x) = 10^{-10} \cdot \I$, the
velocity field $b(x)=(1,-0.5)^T$ and the reaction term $c(x)=1$ to be
constant over the whole domain. The Dirichlet boundary condition is
given by
\begin{equation*}
    g(t,x) = 
    \begin{cases*}
        |\sin{(2t)}|, & \mbox{if } $x \in \{0\}\times\left(\frac{1}{4},\frac{3}{4}+\frac{1}{4}|\sin{(2t)}|\right) $ \\
        0, & \mbox{else}
    \end{cases*}
\end{equation*}
and the initial distribution is given by $u_0(x) := g(0,x) = 0$.
\begin{figure}
    \centering
  \begin{subfigure}[t]{0.48\linewidth}
    \centering
        \includegraphics[height=\linewidth]{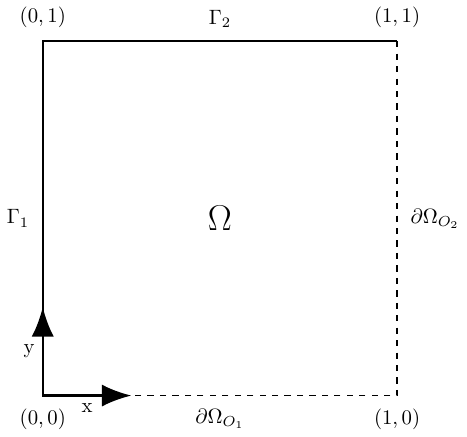}
        \caption{Computational domain of the linear convection-diffusion example. Dirichlet boundary is $\partial \Omega_D = \Gamma_1 \cup \Gamma_2$ and outflow boundary $\partial \Omega_O = \partial \Omega_{O_1} \cup \partial \Omega_{O_2}$.}
        \label{fig:convdiffgeo}
    \end{subfigure}
    \hfill
    \begin{subfigure}[t]{0.48\linewidth}
    \centering
        \includegraphics[height=\linewidth]{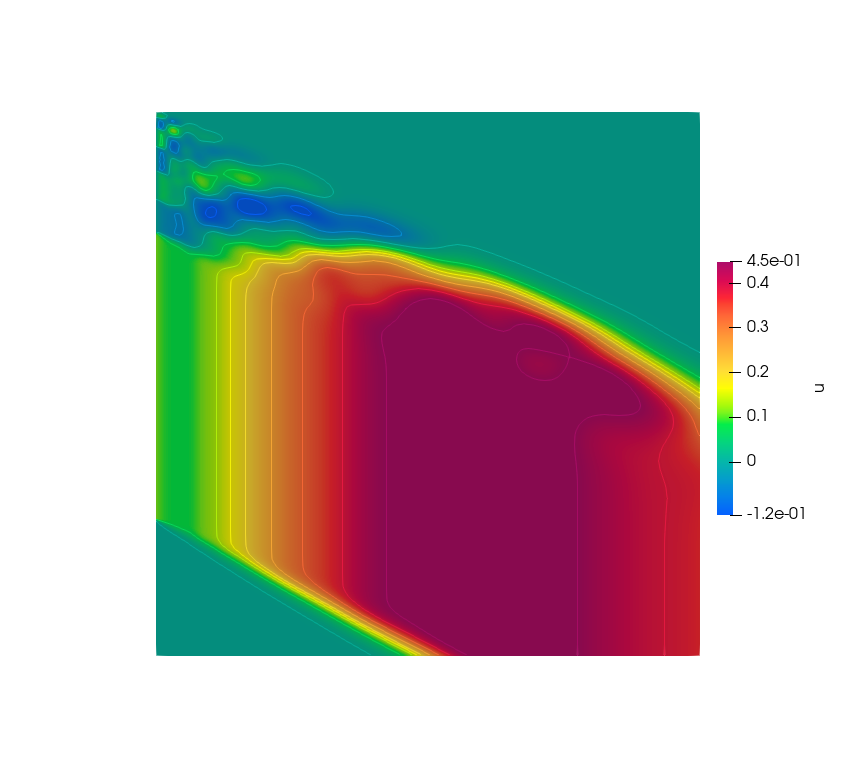}
        \caption{
        Numerical solution for the linear convection-diffusion example at the final time.
    \label{fig:LinearCCFVResults}
    }
    \end{subfigure}
    \caption{Geometry (left) and numerical solution (right) of the linear convection-diffusion example.}
\end{figure}

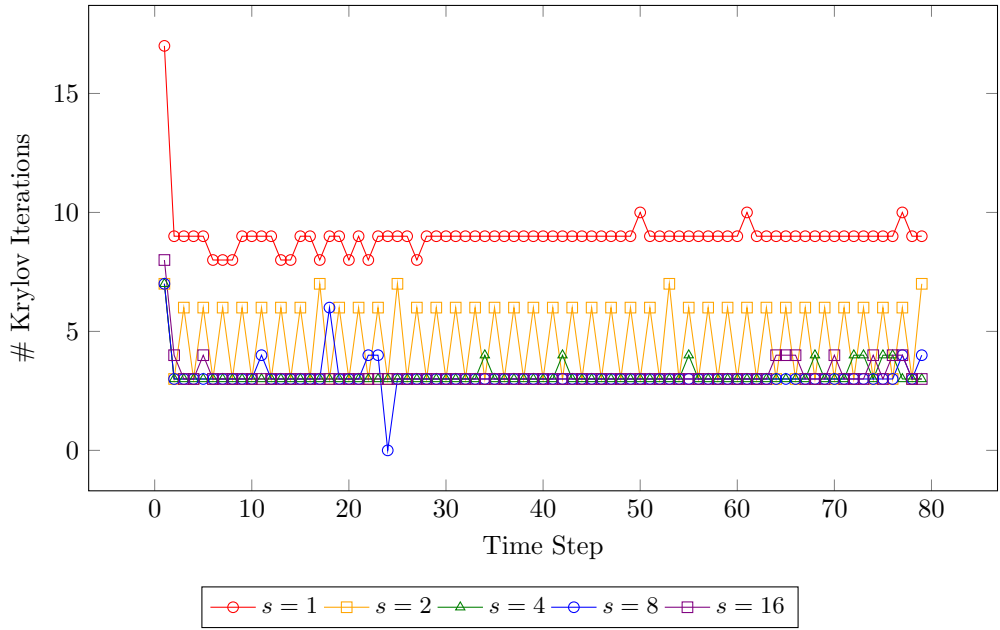
\begin{figure}[htp]
  \centering
    \begin{tikzpicture}
        \begin{axis}
            [
                width=0.85\linewidth,
                height=0.5\linewidth,
                legend to name = timeStepLegend,
                legend columns = 5,
                legend style = {font=\small},
                ylabel = \# Krylov Iterations,
                xlabel = Time Step,
            ]
            \addplot table [x=time_steps, y=ls_iterations, col sep=comma] {convdiffstep1.csv};
            \addlegendentry{$s=1$}

            \addplot table [x=time_steps, y=ls_iterations, col sep=comma] {convdiffstep2.csv};
            \addlegendentry{$s=2$}

            \addplot table [x=time_steps, y=ls_iterations, col sep=comma] {convdiffstep4.csv};
            \addlegendentry{$s=4$}

            \addplot table [x=time_steps, y=ls_iterations, col sep=comma] {convdiffstep8.csv};
            \addlegendentry{$s=8$}

            \addplot table [x=time_steps, y=ls_iterations, col sep=comma] {convdiffstep16.csv};
            \addlegendentry{$s=16$}

        \end{axis}
      \end{tikzpicture}

      \medskip
      \pgfplotslegendfromname{timeStepLegend}

      \caption{The total amount of inner Krylov iterations needed for each time step when calculating $s$ consecutive time steps in parallel. \label{fig:KrylovIterLinearCase}
        }
\end{figure}

For the spatial discretization we used a cell-centered Finite Volume
method (CCFV) on a structured rectangular grid with $640 \times 640$
elements.
The Crank-Nicolson method is used to discretize in time with a time step width of $\tau=0.12$. Simulations were run till the final time $T=4.8$.
The solution at the final time generated by our simulations is shown in Figure~\ref{fig:LinearCCFVResults}.
If only a single time step is used, the inner and outer tolerances
are set to $1e-8$ which means that every time step converges after
a single outer iteration as one would expect. The same could be
done when multiple time steps are considered. But in that case we
do not expect a huge benefit from our approach since the
right-hand-side of the system will be incorrect for all following
time steps. Hence, it should be updated frequently. Currently, the
only possible way to update the right-hand-side of the system is
by adding outer iterations. Therefore, we decided to reduce the
inner tolerance $1e-2$ while the outer tolerance remains
untouched. This means that additional outer iterations are
performed for the first time step that allow to update the
right-hand-side of the following frequently and result in better
preconditioned initial guesses. Another option would be to apply
the update of the right-hand-side directly in the inner
Krylov-solver, which is part of current research. 
The matrix equation arising in the outer iteration is solved using a vectorized
\textsc{Block-BiCGStab} implementation from \textsc{DUNE-ISTL} with an 
AMG preconditioner. Table~\ref{table:ResultsLinear} shows the total number of Newton and Krylov iterations needed 
and the measured solver times. Figure~\ref{fig:KrylovIterLinearCase} 
shows the number of Krylov iterations needed per time depending on the time parallelization.

\pgfplotstableread[
    col sep=comma,
    ]
    {convectiondiffusionresults.csv}\linearresults
\pgfplotstablesort[
    sort key={s},
    sort cmp={int <}]
    \linearresults\linearresults

\begin{table}[htp]
    \centering
        \pgfplotstabletypeset
        [
            col sep=comma,
            columns ={s,nl_iterations,ls_iterations,ls_time,prec_setup_time},
            columns/s/.style={column name=\textbf{$s$}},
            columns/nl_iterations/.style={column name=\bfseries Newton It.},
            columns/ls_iterations/.style={column name=\bfseries Krylov It.},
            columns/ls_time/.style={
                column name=\bfseries Solver Time,
                postproc cell content/.append style={@cell content/.add={}{~[s]}},
                precision=3},
            columns/prec_setup_time/.style={
                column name=\bfseries Prec. Setup Time,
                postproc cell content/.append style={@cell content/.add={}{~[s]}},
                precision=3},
            every head row/.style={before row=\toprule, after row=\midrule},
            every last row/.style={after row=\bottomrule},
        ]
        {\linearresults}
    \centering
    \caption{Results for the linear problem for $s$ simultaneously 
    treated time steps of Crank-Nicolson with a total number of $409600$ degrees of freedom. 
    All results were computed sequentially on 20 cores on a single socket. 
    Computation times are displayed in seconds. }
    \label{table:ResultsLinear} 
\end{table}

\FloatBarrier

\subsection{Nonlinear Diffusion-Reaction}

\begin{figure}
    \centering
    \includegraphics[width=0.8\linewidth
        ]{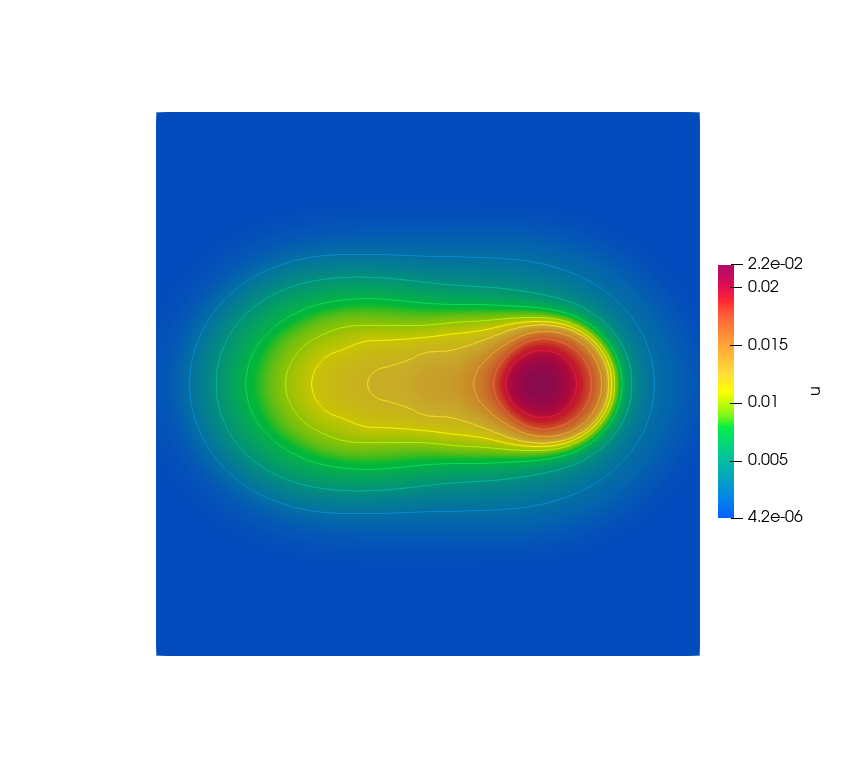}
    \caption{Numerical solution for the nonlinear diffusion-reaction 
    example at the final time.
    \label{fig:NonlinearDiffusionReactionResults}}
\end{figure}

The second example is a nonlinear diffusion-reaction equation 
with Neumann boundary conditions on the unit square $\Omega = [0,1]^2$:
\begin{align*}
    \partial_t u -  \nabla\cdot((1-u)u\nabla u) + \beta (1-u)u &=f &\mbox{in } \Omega \times (0,T], \\
    \partial_n u(t,x) & = 0 & \mbox{on } \partial\Omega \times (0,T],\\
    u(0,x) & = u_0(x) & \mbox{in } \Omega,
\end{align*}
the initial condition
is given by $u_0(x) = 0.5$. The parameter $\beta$ is set to
$1$ and the source term $f$ was given by a moving circle
\begin{equation*}
    f(t,x) = \begin{cases*}
        -0.1, & \mbox{if } $\| x - c \|_2 \leq 0.1$ \\
        0, & \mbox{else},
    \end{cases*}
\end{equation*}
where the center of the moving circle is given by
$c=(0.5 + 0.25 \cos{(6t)} , 0.5)^T$. For the spatial discretization we
chose $Q_1$-Finite Elements on a $640\times640$ rectangular grid.
For the temporal discretization Crank-Nicolson with a time step size of $\tau=0.12$ is used and the simulation was run till $T=9.6$.
The tolerance for the outer iterations was set to $1e-8$ and the inner tolerance to $1e-5$. In contrast to the linear case we chose the same inner tolerance for all experiments since due to the nonlinearity we do not need to enforce more outer iterations to update the right-hand-side of the inner system as we did in the linear case.
The numerical solution at the final time is shown in
Figure~\ref{fig:NonlinearDiffusionReactionResults}. 
Table~\ref{table:ResultsNonLinear} shows the total number of Newton 
and Krylov iterations needed for the whole simulation as well as the 
measured solver times. All linear systems were solved using the 
\textsc{Block-CG} implementation from \textsc{Dune-ISTL} with an AMG 
preconditioner. In Figure~\ref{fig:IterationsDiffusionReaction}
both the inner and outer iterations needed for every time step are displayed.

\pgfplotstableread[
    col sep=comma,
    ]    {nonlineardiffusionreactionresults.csv}\nonlinearresults
\pgfplotstablesort[
    sort key={s},
    sort cmp={int <}]
    \nonlinearresults\nonlinearresults

\begin{table}
  \centering
        \pgfplotstabletypeset
        [
            col sep=comma,
            columns ={s,nl_iterations,ls_iterations,ls_time,prec_setup_time},
            columns/s/.style={column name=\textbf{$s$}},
            columns/nl_iterations/.style={column name=\bfseries Newton It.},
            columns/ls_iterations/.style={column name=\bfseries Krylov It.},
            columns/ls_time/.style={
                column name=\bfseries Solver Time,
                postproc cell content/.append style={@cell content/.add={}{~[s]}},
                precision=3},
            columns/prec_setup_time/.style={
                column name=\bfseries Prec. Setup Time,
                postproc cell content/.append style={@cell content/.add={}{~[s]}},
                precision=3},
            every head row/.style={before row=\toprule, after row=\midrule},
            every last row/.style={after row=\bottomrule},
        ]
        {\nonlinearresults}
    \centering
    \caption{Results for the nonlinear diffusion-reaction equation for 
    $s$ simultaneously treated time steps on a mesh with $640\times640$ 
    elements. All results were computed sequentially on 20 cores on a 
    single socket. Computation times are displayed in seconds.}
    \label{table:ResultsNonLinear} 
\end{table}

\begin{figure}[htp]
    \centering
    \begin{subfigure}[t]{0.85\textwidth}
        \centering
        \begin{tikzpicture}
            \begin{axis}
                [
                    height=0.4\linewidth,
                    legend to name = sharedLegend,
                    legend columns = 5,
                    legend style = {font=\small},
                    ylabel = \# Newton Iterations,
                    xlabel = Time Step,
                    yticklabel style = {
                        text width = 2em,
                        align = right,
                    },
                ]
    
                \addplot table [x=time_steps, y=nl_iterations, col sep=comma] {diffreactstep1.csv};
                \addlegendentry{$s=1$}

                \addplot table [x=time_steps, y=nl_iterations, col sep=comma] {diffreactstep2.csv};
                \addlegendentry{$s=2$}

                \addplot table [x=time_steps, y=nl_iterations, col sep=comma] {diffreactstep4.csv};
                \addlegendentry{$s=4$}

                \addplot table [x=time_steps, y=nl_iterations, col sep=comma] {diffreactstep8.csv};
                \addlegendentry{$s=8$}

                \addplot table [x=time_steps, y=nl_iterations, col sep=comma] {diffreactstep16.csv};
                \addlegendentry{$s=16$}
            \end{axis}
        \end{tikzpicture}
        \caption{Newton iterations}
        \label{fig:NewtonItDiffusionReaction}
    \end{subfigure}
    
    \medskip

    \begin{subfigure}[t]{0.85\textwidth}
        \centering
        \begin{tikzpicture}
            \begin{axis}
                [
                    height=0.4\linewidth,
                    legend style = {draw=none},
                    ylabel = \# Krylov Iterations,
                    xlabel = Time Step,
                    yticklabel style = {
                        text width = 2em,
                        align = right,
                    },
                ]    
                \addplot table [x=time_steps, y=ls_iterations, col sep=comma] {diffreactstep1.csv};
                \addplot table [x=time_steps, y=ls_iterations, col sep=comma] {diffreactstep2.csv};
                \addplot table [x=time_steps, y=ls_iterations, col sep=comma] {diffreactstep4.csv};
                \addplot table [x=time_steps, y=ls_iterations, col sep=comma] {diffreactstep8.csv};
                \addplot table [x=time_steps, y=ls_iterations, col sep=comma] {diffreactstep16.csv};
            \end{axis}
        \end{tikzpicture}
        \caption{Krylov iterations}
        \label{fig:KrylovItDiffusionReaction}    
      \end{subfigure}

      \medskip
      \pgfplotslegendfromname{sharedLegend}

    \caption{Change of iterations (Nonlinear \& Krylov) over the
      different time steps of the nonlinear diffusion-reaction problem, when considering $s$ consecutive time steps in parallel.}
    \label{fig:IterationsDiffusionReaction}
\end{figure}

\FloatBarrier

\subsection{Richards' equation}\label{subsec:Richards}

\begin{figure}[htp]
    \centering
    \begin{subfigure}[t]{0.45\textwidth}
        \centering
        \includegraphics[width=\linewidth,height=\linewidth]{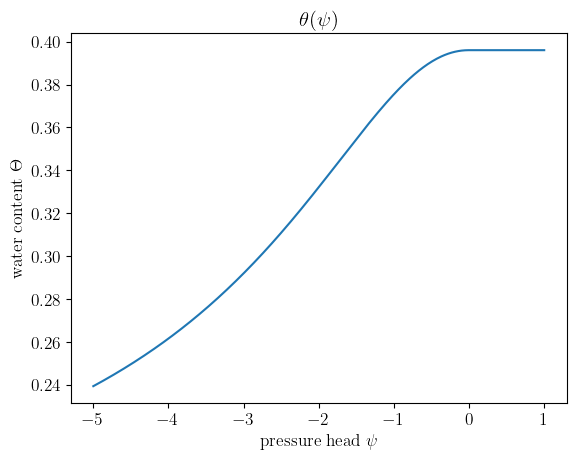}
    \end{subfigure}
    ~
    \begin{subfigure}[t]{0.45\textwidth}
        \centering
        \includegraphics[width=\linewidth,height=\linewidth]{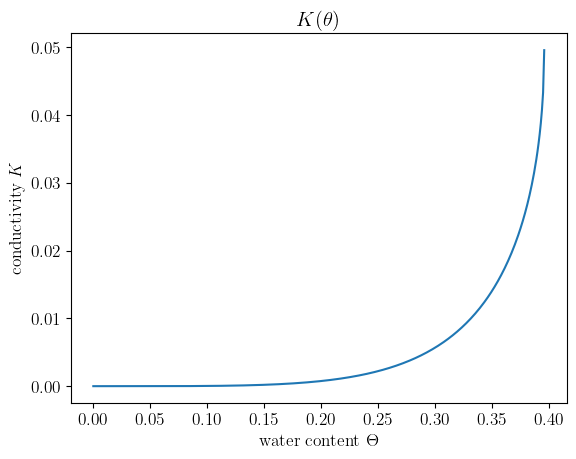}
    \end{subfigure}
    \caption{Relationships between $K,\Theta,\psi$ based on the van 
    Genuchten-Mualem model~\eqref{eq:RelationshipsThetaK}. Soil parameters 
    are taken from Table~\ref{table:Paramters for the Benchmark problem}.}
    \label{fig:WaterContent-Conductivity}
\end{figure}

The last example is Richards' equation. It models unsaturated
flow in soil and reads as
\begin{align*}
    \begin{aligned}
    \partial_t \Theta(\psi) - \nabla \cdot (K(\psi) \nabla (\psi + z)) & = f & \mbox{in } \Omega \times (0,T],\\
    \psi & = g(t,\cdot) & \mbox{on } \Gamma_D \times (0,T],\\
    - K (\psi )(\nabla\psi + e_z) \cdot n & = 0 & \mbox{on } \Gamma_N \times (0,T] \\
    \psi(0,\cdot) & = \psi_0(\cdot) & \mbox{in } \Omega,
    \end{aligned}
\end{align*}
where the primal variable $\psi$ denotes the pressure head, $\Theta(\psi)$ 
the water content, $K(\psi)$ the conductivity, $z$ the height against 
gravitation, $e_z:=\nabla z$ and $f$ is a source term. The highly nonlinear 
dependencies between $K,\Theta$ and $\psi$ are described based on 
experimental results, and we use the widely used model of 
\cite{van1980closed}

\begin{equation}
    \begin{aligned}
        \theta(\psi) &= 
        \begin{cases}
            \theta_R + (\theta_S-\theta_R)\left(\frac{1}{1+(-\alpha\psi)^n}
            \right)^{\frac{n-1}{n}} & \text{for } \psi \leq 0\\
            \theta_S & \text{for } \psi > 0
        \end{cases}\\
        K(\psi) &= K_S\left(\frac{\theta(\psi)}{\theta_S}\right)^{\frac12}
        \left[1-\left(1-\left(\frac{\theta(\psi)}{\theta_S}\right)^{\frac{n}{n-1}}\right)^{\frac{n-1}{n}}\right]^2,
    \end{aligned}
    \label{eq:RelationshipsThetaK}
\end{equation}
where $\Theta_S$ denotes the water content and $K_S$ the 
hydraulic conductivity in the fully saturated medium. $\Theta_R$ 
stands for the residual water content and $\alpha$ and $n$ are model 
parameters related to the soil properties. The relationships are shown 
in Figure~\ref{fig:WaterContent-Conductivity}.

\begin{figure}
    \centering
  \begin{subfigure}[t]{0.48\linewidth}
    \centering
        \includegraphics[height=4.5cm]{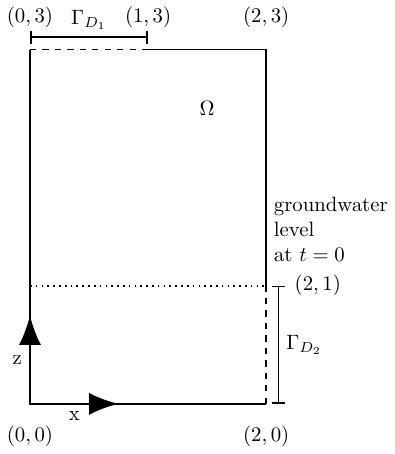}
        \caption[Geometry for the silt loam benchmark problem]{Geometry for the silt loam benchmark problem, from \cite{list2016study}.}
        \label{fig:benchmark-problem-geometry}
    \end{subfigure}
    \hfill
    \begin{subfigure}[t]{0.48\linewidth}
    \centering
        \includegraphics[height=4.5cm]{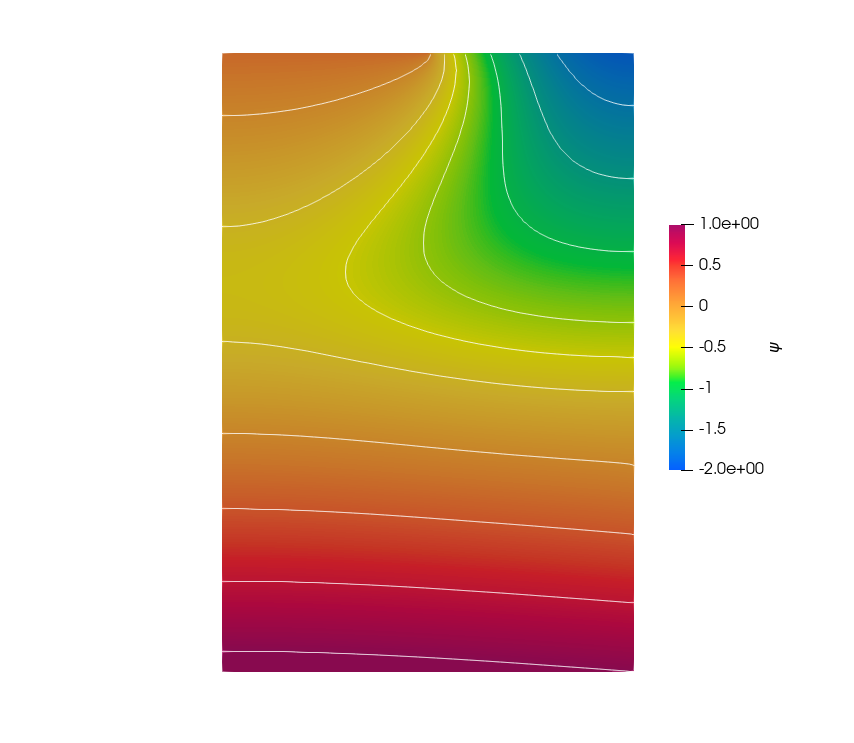}
        \caption{Numerical solution for the pressure height $\psi$ after $4.5h$.}
        \label{fig:richardsfinaltime}
    \end{subfigure}
    \caption{Geometry (left) and numerical solution at the final time (right) for the benchmark problem of the Richards' equation.}
\end{figure}

In the considered example the Richards' equation will be nonlinear and 
degenerate elliptic-parabolic. Its solutions typically will have low 
regularity \parencite{alt1983quasilinear}. The exact setting we use here is 
a benchmark problem used by e.g. \cite{schneid2000hybrid,list2016study} 
amongst others. It models the recharge of a groundwater reservoir from 
a drainage trench in two spatial dimensions. The domain $\Omega \subset \mathbb{R}^2$ 
represents a vertical cross-section of the subsurface. The groundwater 
table is fixed by a Dirichlet boundary condition for the pressure height 
for $z \in [0,1]$ on the right-hand-side of $\Omega$, $\Gamma_{D_2}$. 
On the upper boundary of $\Omega$ $\Gamma_{D_1}$ for $x \in [0,1]$ the 
drainage trench is modelled by a transient Dirichlet condition. On the 
remaining part of the boundary, no-flow conditions are imposed. The left 
part of the domain boundary can be interpreted as symmetry axis of the 
geometry and the lower part of the domain boundary can be seen as transition 
to an aquitard. The exact geometry is given by

\begin{align*}
    \Omega & = (0,2) \times (0,3), \\
    \Gamma_{D_1} & = \left\lbrace (x,z) \in \partial \Omega \, | \, x \in [0,1] \wedge z = 3 \right\rbrace, \\
    \Gamma_{D_2} & = \left\lbrace (x,z) \in \partial \Omega \, | \, x = 2 \wedge z \in [0,1] \right\rbrace , \\
    \Gamma_D & = \Gamma_{D_1} \cup \Gamma_{D_2}, \\
    \Gamma_N & = \partial \Omega \setminus \Gamma_D,
\end{align*}
and is illustrated in Figure~\ref{fig:benchmark-problem-geometry}.

The initial and boundary conditions are given by
\begin{align*}
    \begin{split}
        \psi ( t,x,z ) = 
        \begin{cases}
            -2 + 2.2 t / \Delta t_D, & \text{on } \Gamma_{D_1} , t \leq \Delta t_D, \\
            0.2, & \text{on } \Gamma_{D_1} , t > \Delta t_D, \\
            1-z, & \text{on } \Gamma_{D_2}, 
        \end{cases} \\
        - K ( \psi ( t,x,z) ) ( \nabla \psi (t,x,z) + e_z ) \cdot n = 0 \quad \text{on } \Gamma_N, \\
        \psi^0 (x,z) = 1 - z \quad \text{on } \Omega.
    \end{split}
\end{align*}

The set of parameters characterizing the soil are taken from 
\cite{van1980closed}. They represent silt loam and are shown in Table~\ref{table:Paramters for the Benchmark problem}.
For the spatial discretization we choose first order Finite Elements on 
a mesh with $400\times600$ elements. For the temporal discretization 
we choose the implicit Euler method. Since Newton's method might fail 
as shown in \cite{radu2006newton}, we use the \emph{L-Scheme} as the 
linearization method as done in \cite{list2016study}. The L-Scheme can 
be interpreted as a quasi-Newton method that uses the monotonicity of 
$\Theta$. It replaces its derivative by a constant $L$ and was proposed 
for the Richards' equation by \cite{slodicka2002robust,pop2004mixed,yong1996numerical}. 
We use a time step size of $\tau=\frac{1}{96}$ and perform the simulations 
until the final time $T=\frac{3}{16}$. The parameter for the time evolution 
of the upper Dirichlet boundary is given by $\Delta t_D = \frac{1}{16}$. 
The L-Scheme parameter $L$ is chosen to be the Lipschitz constant of 
$\Theta$, that is $L=4.501\cdot10^{-2}$. In this example the time unit 
is $1$ day and all spatial dimensions are given in meters.
For the stopping criterion of the outer L-Scheme iteration we again
refer to \cite{list2016study}, where a combination of absolute and relative
error is used with absolute and relative tolerances of $1e-5$.

\begin{table}[htp]
    \centering
    \begin{tabular}{ll}
        \toprule 
        & \bfseries Silt loam benchmark problem\\ 
        \midrule
        \bfseries van Genuchten parameters:  &   \\ 
        $\Theta_S$ & $0.396$  \\ 
        $\Theta_R$ & $0.131$ \\ 
        $\alpha$ & $0.423$ \\ 
        $n$ & $2.06$ \\ 
        $K_S$ & $4.96 \cdot 10^{-2}$ \\
        \bottomrule
    \end{tabular} 
    \caption[Parameters of the example problems]{Soil parameters for the silt 
    loam benchmark problem \cite{van1980closed}}\label{table:Paramters for the Benchmark problem}
\end{table}

\begin{figure}
    \centering
    \begin{subfigure}[t]{0.85\textwidth}
        \begin{tikzpicture}
            \begin{axis}
                [
                    height=0.4\linewidth,
                    legend to name = sharedLegend,
                    legend columns = 5,
                    legend style = {font=\small},
                    ylabel = \# L-Scheme Iterations,
                    xlabel = Time Step,
                    yticklabel style = {
                        text width = 2em,
                        align = right,
                    },
                ]
    
                \addplot table [x=time_steps, y=nl_iterations, col sep=comma] {richards1.csv};
                \addlegendentry{$s=1$}

                \addplot table [x=time_steps, y=nl_iterations, col sep=comma] {richards2.csv};
                \addlegendentry{$s=2$}

                \addplot table [x=time_steps, y=nl_iterations, col sep=comma] {richards4.csv};
                \addlegendentry{$s=4$}

                \addplot table [x=time_steps, y=nl_iterations, col sep=comma] {richards8.csv};
                \addlegendentry{$s=8$}

                \addplot table [x=time_steps, y=nl_iterations, col sep=comma] {richards16.csv};
                \addlegendentry{$s=16$}
            \end{axis}
        \end{tikzpicture}
        \caption{L-Scheme iterations}
        \label{fig:L-SchemeIterationsRichardsBenchmark}
    \end{subfigure}
    \bigskip
    \begin{subfigure}[t]{0.85\textwidth}
        \begin{tikzpicture}
            \begin{axis}
                [
                    height=0.4\linewidth,
                    legend style = {draw=none},
                    ylabel = \# Krylov Iterations,
                    xlabel = Time Step,
                    yticklabel style = {
                        text width = 2em,
                        align = right,
                    },
                ]
    
                \addplot table [x=time_steps, y=ls_iterations, col sep=comma] {richards1.csv};
                \addplot table [x=time_steps, y=ls_iterations, col sep=comma] {richards2.csv};
                \addplot table [x=time_steps, y=ls_iterations, col sep=comma] {richards4.csv};
                \addplot table [x=time_steps, y=ls_iterations, col sep=comma] {richards8.csv};
                \addplot table [x=time_steps, y=ls_iterations, col sep=comma] {richards16.csv};
            \end{axis}
        \end{tikzpicture}
        \caption{Krylov iterations}
        \label{fig:KrylovIterationsRichardsBenchmark}    
      \end{subfigure}

      \medskip
      \pgfplotslegendfromname{sharedLegend}

    \caption{Change of iterations (L-Scheme \& Krylov) over the
      different time steps of the Richards problem, when considering $s$ 
      consecutive time steps in parallel.}
    \label{fig:LSchemeKrylovIterationsRichards}
\end{figure}
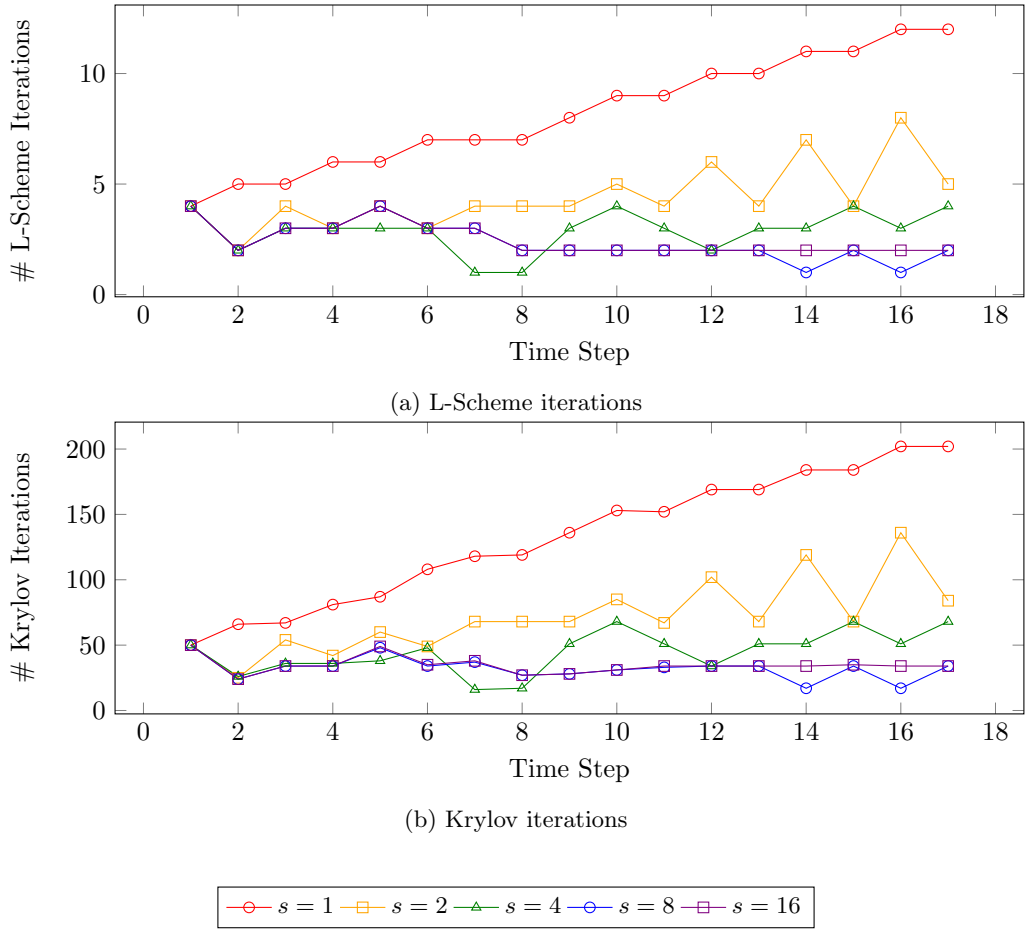

Figure~\ref{fig:richardsfinaltime} shows the numerical solution at 
the final time computed with our method and Table~\ref{table:ResultsRichards} 
shows the total number of L-Scheme iterations and inner Krylov iterations 
needed together with the total time needed by the linear solver. All 
arising linear systems were solved using a vectorized \textsc{Block-CG} 
implementation with an AMG preconditioner.
Figure~\ref{fig:LSchemeKrylovIterationsRichards} shows the evolution of 
the L-Scheme iterations and inner Krylov iterations needed for every time step.

\pgfplotstableread[
    col sep=comma,
    ]
    {richardsresults.csv}\richardsresults
\pgfplotstablesort[
    sort key={s},
    sort cmp={int <}]
    \richardsresults\richardsresults

\begin{table}[htp]
  \centering
        \pgfplotstabletypeset
        [
            col sep=comma,
            columns ={s,nl_iterations,ls_iterations,ls_time,prec_setup_time},
            columns/s/.style={column name=\textbf{$s$}},
            columns/nl_iterations/.style={column name=\bfseries L-Scheme It.},
            columns/ls_iterations/.style={column name=\bfseries Krylov It.},
            columns/ls_time/.style={
                column name=\bfseries Solver Time,
                postproc cell content/.append style={@cell content/.add={}{~[s]}},
                precision=3},
            columns/prec_setup_time/.style={
                column name=\bfseries Prec. Setup Time,
                postproc cell content/.append style={@cell content/.add={}{~[s]}},
                precision=3},
            every head row/.style={before row=\toprule, after row=\midrule},
            every last row/.style={after row=\bottomrule},
        ]
        {\richardsresults}
    \centering
    \caption{Results for the benchmark problem of the Richards' equation for 
    $s$ simultaneously treated time steps on a mesh with $400\times600$ 
    elements. All results were computed sequentially on 20 cores on a 
    single socket. Computation times are displayed in seconds.}
    \label{table:ResultsRichards} 
\end{table}

\FloatBarrier

\subsection{Discussion of the results}
The presented numerical experiments show that the total number of linear Krylov 
iterations and, in the nonlinear case, also the number of Newton 
iterations has been significantly reduced when computing multiple 
consecutive time steps in parallel. In the linear case
the number of linear iterations was reduced by approximately $50-60\%$
(see Table~\ref{table:ResultsLinear}). As mentioned earlier, we chose to solve the inner system less accurate when considering multiple time steps in parallel in the linear case. This results in more outer iterations needed for the first time step and consequently the right-hand-side is updated more frequently resulting in a better preconditioned next time step. This can especially be seen for $s=2$ where the number of outer iterations increases and at the same time the number of inner iterations decreases significantly which yields to an overall speedup of $1.8$ which is nearly optimal.

In the nonlinear case one observes that the total number of inner and outer iterations were reduced with more parallelization, i.e. for the diffusion-reaction example both iteration counts were overall reduced by $60-70\%$. For $s=2$ the obtained reduction of $48\%$ is nearly optimal (see Table~\ref{table:ResultsNonLinear}). Larger $s$ only yield a moderate benefit which can be explained by Figure~\ref{fig:NewtonItDiffusionReaction} where it is shown that even for $s=2$ some time steps already converge after a single outer iteration and hence there cannot be much more benefit.
Similar results can be seen for Richards' equation (Table~\ref{table:ResultsRichards}). Again outer and inner iteration counts were reduced by $60-70\%$. As before, for $s=2$ the reduction in inner and outer iterations is nearly optimal but for larger $s$ the benefit is only moderate. Most likely this is caused by the larger splitting error and by the influence of the nonlinearity.

Figures~\labelcref{fig:KrylovIterLinearCase,fig:IterationsDiffusionReaction,fig:LSchemeKrylovIterationsRichards} show the evolution of inner and in the nonlinear case outer iterations per time step. After populating the pipeline, all experiments show a reduction in both iteration counts. The reason for this is that a time step, when newly shifted to the first column during the pipelining procedure, has an improved initial guess due to the previous approximate Newton updates before. However, this effect depends on the nonlinearity which can be seen e.g. for the Richards' equation. Although, for larger $s$ we still get a reduction in the iteration numbers, the benefit is only moderate.

Figure~\ref{fig:NewtonItDiffusionReaction} shows that 
in the nonlinear diffusion-reaction 
example one actually reaches the behaviour mentioned in
Remark~\ref{remark:non-linear-convergence}, with sufficient
parallelization the pipelining procedure only requires one Newton iteration per timestep.

Looking at the linear solver times in Tables~\labelcref{table:ResultsLinear,table:ResultsNonLinear,table:ResultsRichards} one observes
that for the nonlinear examples the best results are obtained for $s=4$. For instance, the increased arithmetic intensity and vectorization yield a speedup of $2.6$ for the Richards' equation and for the 
preconditioner set up time a speedup of $2.9$ is reached. In the linear case the best speedup is gained for $s=2$.
The optimal speedup is limited by the maximum of $8$ vector lanes on \emph{Intel Xeon Gold 6148} with \emph{AVX512}. For values of $s$ larger than $4$ the measured times increase and in the nonlinear diffusion-reaction example for $s=16$ it is even more time needed than for $s=1$, although the least iterations are needed. This is supposed to be caused by the worse cache behaviour mentioned in Section~\ref{sec:PerfomanceModelling}.

	\section{Conclusion}\label{sec:Conclusion}
In this paper we presented a new generic framework for stiffly accurate diagonally implicit
Runge-Kutta methods based on Parallel in Time ideas that makes use of 
the advantages of Block Krylov solvers. Furthermore, multiple properties 
of modern hardware like the use of SIMD-instructions and the efficient 
use of memory-bandwidth are respected. This is done by transforming a 
system for multiple consecutive time steps or Runge-Kutta stages into 
a single matrix equation that can be solved using a vectorized linear 
solver. By this the arithmetic intensity in the solver was increased. 
Hence, the framework is particularly well suited for low-order methods 
that are widely used in applications to porous media problems.
The results show that the solver time of the linear solver was reduced by a factor of $2-3$.
The {total} number of linear solver iterations was {reduced} by around $60-70\%$. 
This means that the communication overhead on distributed systems would
also be reduced by this factor. Due to its structure our approach is 
applicable for a wide class of linear and nonlinear problems.

Two aspects are important regarding the performance of the
  method and thus its applicability to other problems.
The first aspect is the splitting error. Reducing it would
  reduce the number of outer iterations
and by
this improve the performance of the method. In particular in the
nonlinear case this question is tightly related to the dynamics of
the nonlinear operator. An analysis of the splitting error in this
case would allow for a better a-priori estimate of the potential
benefit of this type of parallelization.
A second aspect is the inherently sequential nature of
  problem, which
becomes
visible in the numerical experiments.
With increasing parallelism, i.e. more time steps, the method becomes effectively slower again.
In PinT methods a common
approach is to use a discretization with coarser time steps and by this
improve the overall convergence. Such ideas are employed in different
ways in most parallel in time methods, c.f.~\parencite{gander2025time}, a
future aim should thus be to also introduce such concepts into the
presented method, in order to allow for more parallel stages.

	\printbibliography

	\section*{Statements and declarations}

	\paragraph{Funding}
	This research was supported by the Deutsche Forschungsgemeinschaft (DFG, German Research Foundation) under Germany's Excellence Strategy EXC 2044 -- 390685587, Mathematics Münster: Dynamics-Geometry-Structure and DFG projects 504505951 (BlockXT - Block methods to transparently accelerate and vectorise time dependent simulations) and 230658507 (EXA-DUNE - Flexible PDE Solvers, Numerical Methods, and Applications).

	\paragraph{Competing interests}
	The authors have no relevant financial or non-financial interests to disclose.

\end{document}